\newif\ifusual

\usualfalse

\ifusual

\documentclass[11pt
 ]{article}

\evensidemargin\oddsidemargin

\else

\documentclass[11pt,twoside,reqno%
 ]{amsart}

\evensidemargin\oddsidemargin

\fi

\usepackage{JPRM}
\usepackage{DLdef}
\theoremstyle{thm}
\newtheorem{Conj}{\secno Conjecture}

\usepackage{soul}
\usepackage{epigraph}
\begin{document}

\ifusual \else
\FirstPageHead{3}{2}{2007}{\pageref{firstpage}--\pageref{lastpage}}

\fi

\thispagestyle{empty}

\Name{Towards classification of simple finite dimensional modular
Lie~superalgebras}

\label{firstpage}

\Author{Dimitry Leites}
\thanks{I am thankful to S.~Bouarroudj, P.~Grozman, A.~Lebedev, I.~Shchepochkina,
and also to V.~Serganova and Yu.~Kochetkov for help; MPIMiS,
Leipzig, for financial support and most creative environment during
2004-06 when I was Sophus-Lie-Professor there.}
\Address{MPIMiS, Inselstr. 22, DE-04103 Leipzig, Germany\\
on leave from Department of Mathematics, University of Stockholm,
Roslagsv. 101, Kr\"aft\-riket hus 6, SE-106 91 Stockholm, Sweden;
mleites@math.su.se}

\begin{abstract} A way to construct (conjecturally all) simple finite
dimensional modular
Lie (super)algebras over algebraically closed fields of
characteristic not 2 is offered. In characteristic 2, the method is
supposed to give only simple Lie (super)algebras graded by integers
and only some of the non-graded ones). The conjecture is backed up
with the latest results computationally most difficult of which are
obtained with the help of Grozman's software package SuperLie.

\keywords {Cartan prolongation, Kostrikin-Shafarevich conjecture,
modular Lie algebra, Lie superalgebra}

\subjclass{17B50, 70F25}
\end{abstract}


\markboth{Dimitry Leites}{On classification of simple Lie
superalgebras}


\epigraph{Characteristic $p$ is for the time when we retire.}
{Sasha Beilinson, when we all were young. }

\vspace*{.5cm}

\section{Introduction}

The purpose of this transcript of the talk to be presented in March
2007 at the 3rd International Conference on 21st Century Mathematics
2007, School of Mathematical Sciences (SMS), Lahore, is to state
problems, digestible to Ph.D. students (in particularly, the
students at SMS) and worth (Ph.D. diplomas) to be studied (without
waiting till retirement time), together even with ideas of their
solution. In the process of formulating the problems, I'll overview
the classical and latest results. To be able to squeeze the material
into the prescribed 10 pages, all background is supplied by
accessible references: We use standard notations of \cite{FH, S};
for a precise definition of the Cartan prolongation and its
generalizations (Cartan-Tanaka-Shchepochkina or CTS-prolongations),
see \cite{Shch}; see also \cite{BGL3}--\cite{BGL5}. Hereafter $\Kee$
is an algebraically closed (unless finite) field, $\Char\Kee=p$.

The works of S.~Lie, Killing and Cartan, now classical, completed
classification  over $\Cee$ of {\bf simple Lie algebras}
\begin{equation}
\label{1} \renewcommand{\arraystretch}{1.0}
\begin{array}{c}
\text{{\bf  of finite dimension} and {\bf certain infinite dimensional} }\\
\text{{\bf  (of polynomial vector fields, or \lq\lq vectorial'' Lie
algebras)}.}
\end{array}
\end{equation}
In addition to the above two types, there are several more
interesting types of simple Lie algebras but they do not contribute
to the solution of our problem: {\bf classification of simple finite
dimensional modular Lie (super)algebras}, except one: the queer type
described below (and, perhaps, examples, for $p=2$, of the types
described in \cite{J, Sh} and their generalizations, if any).
Observe that all finite dimensional simple Lie algebras are of the
form $\fg(A)$; for their definition embracing the modular case and
the classification, see \cite{BGL5}.

Lie algebras and Lie superalgebras over fields in characteristic
$p>0$, a.k.a. {\it modular} Lie (super)algebras, were distinguished
in topology in the 1930s. The {\bf simple} Lie algebras drew
attention (over finite fields $\Kee$) as a byproduct of
classification of simple finite groups, cf. \cite{St}. Lie {\bf
super}algebras, even simple ones, did not draw much attention of
mathematicians until their (outstanding) usefulness was observed by
physicists in the 1970s. Researchers discovered more and more of new
examples of simple modular Lie algebras for decades until Kostrikin
and Shafarevich (\cite{KSh}) formulated a {\bf conjecture} embracing
all previously found examples for $p>7$. The {\bf generalized
KSh--conjecture} states (for a detailed formulation, convenient to
work with, see \cite{Leb2}):

{\sl Select a $\Zee$-form $\fg_\Zee$ of every $\fg$ of
type\footnote{Notice that the modular analog of the polynomial
algebra---the algebra of divided powers---{\it and all prolongs}
(vectorial Lie algebras) acquire for $p>0$ one more (shearing)
parameter $\underline{N}$.} $(\ref{1})$, take
$\fg_\Kee:=\fg_\Zee\otimes_\Zee\Kee$ and its simple subquotient
$\fs\fii(\fg_\Kee)$ (for the Lie algebras of vector fields, there
are several, depending on $\underline{N}$). Together with
deformations\footnote{It is not clear, actually, if the conventional
description of infinitesimal deformations in terms of $h^2(\fg:
\fg)$ can always be applied if $p>0$. This concerns both Lie
algebras and Lie superalgebras (for the arguments, see \cite{LL});
to give the correct (better say, universal) notion is an open
problem, but we let it pass for the moment, besides, for $p\neq 2$
and $\fg$ with Cartan matrix, the conventional interpretation is
applicable, see \cite{BGL4}.} of these examples we get in this way
all simple finite dimensional Lie algebras over algebraically closed
fields if $p>5$. If $p=5$, Melikyan's examples\footnote{For their
description as prolongs, and newly discovered super versions, see
\cite{GL4,BGL3}.} should be added to the examples obtained by the
above method.}

After 30 years of work of several teams of researchers, Block,
Wilson, Premet and Strade proved the generalized KSh conjecture
for $p>3$, see \cite{S}.

Even before the KSh conjecture was proved, its analog was offered in
\cite{KL} for $p=2$. Although the KL conjecture was, as is clear
now, a bit overoptimistic (in plain terms: wrong, as stated), it
suggested a way to get such an abundance of examples (to verify
which of them are really simple is one of the tasks still open) that
Strade \cite{S} cited \cite{KL} as an indication that the case $p=2$
is too far out of reach by modern means\footnote{Contrarywise, {\bf
the \lq\lq punch line" of this talk is: Cartan did not have the
modern root technique, but got the complete list of simple Lie
algebras; let's use his \lq\lq old-fashioned" methods: they work!}
Conjecture 2 expresses our hope in precise terms. How to {\sl prove}
the completeness of the list of examples we will have unearthed is
another story.}. Still, \cite{KL} made two interesting observations:
It pointed at a striking similarity (especially for $p=2$) between
modular Lie algebras and Lie {\bf super}algebras (even over $\Cee$),
and it introduced totally new characters --- {\it Volichenko
algebras} (inhomogeneous with respect to parity subalgebras of Lie
superalgebras); for the classification of simple Volichenko algebras
(finite dimensional and infinite dimensional vectorial) over $\Cee$,
see \cite{LSer} (where one of the most interesting examples is
missed, the version of the proof with repair will be put in
\texttt{arXiv} soon).

Recently Strade had published a monograph \cite{S} summarizing the
description of newly classified simple finite dimensional Lie
algebras over the algebraically closed fields $\Kee$ of
characteristic $p>3$, and also gave an overview of the \lq\lq
mysterious" examples (due to Brown, Frank, Ermolaev and Skryabin) of
simple finite dimensional Lie algebras for $p=3$ with no
counterparts for $p>3$. Several researchers started afresh to work
on the cases where $p=2$ and 3, and new examples of simple Lie
algebras with no counterparts for $p\neq 2, 3$ started to appear
(\cite{J,GL4,Leb1}, observe that the examples of \cite{GG,Lin1} are
erroneous as observed in MathRevies and \cite{Leb2}, respectively).
The \lq\lq mysterious" examples of simple Lie algebras for $p=3$
were interpreted as vectorial Lie algebras preserving certain
distributions (\cite{GL4}).

While writing \cite{GL4} we realized, with considerable dismay, that
there are reasons to put to doubt the universal applicability of the
conventional definitions of the enveloping algebra $U(\fg)$ (and its
restricted version) of a given Lie algebra $\fg$, and hence doubt in
applicability of the conventional definitions of Lie algebra
representations and (co)homology to the modular case, cf. \cite{LL}.
But even accepting conventional definitions, there are plenty of
problems to be solved before one will be able to start writing the
proof of classification of simple modular Lie algebras, to wit:
describe irreducible representations (as for vectorial Lie
superalgebras, see \cite{GLS}), decompose the tensor product of
irreducible
representations into indecomposables, cf. \cite{Cla}, 
and many more; for a review, see \cite{GL2}.

Classification of simple Lie {\bf super}algebras for $p>0$ and the
study of their representations are of independent interest. A
conjectural list of simple finite dimensional Lie superalgebras over
an algebraically closed fields $\Kee$ for $p>5$, known for some
time, was recently cited in \cite{BjL}:

\begin{Conj}[Super KSh, $p>5$]  Apply the steps of the
KSh conjecture to the simple complex Lie {\bf super}algebras $\fg$
of types $(\ref{1})$. The examples thus obtained exhaust all simple
finite dimensional Lie superalgebras over algebraically closed
fields if $p>5$.
\end{Conj}

The examples obtained by this procedure will be referred to as {\it
KSh-type Lie superalgebras}. The first step towards obtaining the
list of KSh-type Lie superalgebras is classification of simple Lie
superalgebras of types $(\ref{1})$ over $\Cee$. This is done by
I.Shchepochkina and me; for summaries with somewhat different
emphases, and proof, see \cite{K2, K3, LSh}.

For a classification of finite dimensional simple modular Lie
algebras with Cartan matrix, see \cite{WK, KWK}. For a
classification of finite dimensional simple modular Lie
superalgebras with Cartan matrix, see \cite{BGL5}. Not all finite
dimensional Lie superalgebras over $\Cee$ are of the form $\fg(A)$;
in addition to them, there are also queer types described below, and
even simple vectorial.

I am sure that the same ideas of Block and Wilson that proved
classification of simple {\it restricted} Lie algebras for $p>5$
will work, if $p>5$, {\it mutatis mutandis}, for Lie superalgebras
and ideas of Premet and Strade will embrace the non-restricted
superalgebras as well; although the definition of {\it
restrictedness} and even of the Lie superalgebra itself acquire more
features, especially for $p=2$.

Here I will describe the cases $p\leq 5$ where the situation is
different and suggest {\it another}, different from KSh, way to get
simple examples.

\section{How to construct simple Lie algebras and superalgebras}

\ssec{How to construct simple Lie algebras if $p=0$} Let us recall
how Cartan used to construct simple $\Zee$-graded Lie algebras
over $\Cee$ of polynomial growth \cite{C} and finite depth. Now
that they are classified (for examples of infinite depth, see
\cite{K}), we know that, {\sl all of them can be endowed with a
$\Zee$-grading $\fg=\mathop{\oplus}\limits_{-d\leq i}\fg_i$ of
depth $d=1$ or 2 so that $\fg_0$ is a simple Lie algebra $\fs$ or
its trivial central extension $\fc\fs=\fs\oplus \fc$, where $\fc$
is a 1-dimensional center. Moreover, simplicity of $\fg$ requires
$\fg_{-1}$ to be an irreducible $\fg_0$-module that generates
$\fg_-:=\mathop{\oplus}\limits_{i<0}\fg_i$ and $[\fg_{-1},
\fg_{1}]=\fg_0$.}

Yamaguchi's theorem \cite{Y}, reproduced in \cite{GL4, BjL}, states
that {\sl for almost all simple finite dimensional Lie algebras
$\fg$ over $\Cee$ and their $\Zee$-gradings
$\fg=\mathop{\oplus}\limits_{-d\leq i}\fg_i$, the generalized Cartan
prolong of $\fg_{\leq}=\mathop{\oplus}\limits_{-d\leq i\leq 0}\fg_i$
is isomorphic to $\fg$, the rare exceptions being two of the four
series of simple vectorial algebras; the other two series being {\bf
partial} prolongs (perhaps, after factorization modulo center).}

For illustration, we construct simple Lie algebras of type
$(\ref{1})$ by induction:

\footnotesize

\underline{{\bf Depth $d=1$}}. Here we use either usual or partial
Cartan prolongations.

1) we start with 1-dimensional $\fc$, so $\dim \fg_{-1}=1$ due to
irreducibility. The complete prolong is isomorphic to $\fvect(1)$,
the partial one to $\fsl(2)$.

2) Take $\fg_0=\fc\fsl(2)=\fgl(2)$ and its irreducible module
$\fg_{-1}$. The component  $\fg_{1}$ of the  Cartan prolong is
nontrivial only if $\fg_{-1}$ is $R(\varphi_1)$ or $R(2\varphi_1)$,
where $\varphi_i$ is the $i$th fundamental weight of the simple Lie
algebra $\fg$ and $R(w)$ is the irreducible representation with
highest weight $w$.

2a) If $\fg_{-1}$ is $R(\varphi_1)$, the component $\fg_1$ consists
of two irreducible submodules, say $\fg_1'$ or $\fg_1''$. We can
take any one of them or both; together with $\fg_{-1}\oplus\fg_0$
this generates $\fsl(3)$ or $\fsvect(2)\subplus \fd$, where $\fd$ is
spanned by an outer derivation, or $\fvect(2)$, respectively.

2b) If $\fgl(2)\simeq\fc\fo(3)\simeq\fc\fsp(2)$-module $\fg_{-1}$
is $R(2\varphi_1)$, then $(\fg_{-1}, \fg_0)_*\simeq\fo(5)\simeq
\fsp(4)$.

3) Induction: Take $\fg_0=\fc\fsl(n)=\fgl(n)$ and its irreducible
module $\fg_{-1}$. The component  $\fg_{1}$ of the Cartan prolong is
nontrivial only if $\fg_{-1}$ is $R(\varphi_1)$ or $R(2\varphi_1)$
or $R(\varphi_2)$.

3a) If $\fg_{-1}=R(\varphi_1)$, then $\fg_1$ consists of two
irreducible submodules, $\fg_1'$ or $\fg_1''$. Take any of them or
both; together with $\fg_{-1}\oplus\fg_0$ this generates $\fsl(n+1)$
or $\fsvect(n)\subplus \fd$, where $\fd$ is spanned by an outer
derivation, or $\fvect(n)$, respectively.

3b) If $\fg_{-1}=R(2\varphi_1)$, then $(\fg_{-1},
\fg_0)_*\simeq\fsp(2n)$.

3c) If $\fg_{-1}=R(\varphi_2)$, then $(\fg_{-1},
\fg_0)_*\simeq\fo(2n)$.

4) The induction with $\fg_0=\fc\fo(2n-1)$-module $R(\varphi_1)$
returns $(\fg_{-1}, \fg_0)_*\simeq\fo(2n+1)$. Observe that
$\fsl(4)\simeq \fo(6)$. The induction with
$\fg_0=\fc\fo(2n)$-module $R(\varphi_1)$ returns $(\fg_{-1},
\fg_0)_*\simeq\fo(2n+2)$. (We have obtained $\fo(2n)$ twice;
analogously, there many ways to obtain other simple Lie algebras
as prolongs.)

5) The $\fg_0=\fsp(2n)$-module $\fg_{-1}=R(\varphi_1)$ yields the
Lie algebra $\fh(2n)$ of Hamiltonian vector fields.

\underline{$\fe(6)$, $\fe(7)$}. The $\fg_0=\fc\fo(10)$-module
$\fg_{-1}=R(\varphi_1)$ yields $\fe(6)$; the
$\fg_0=\fc\fe(6)$-module $\fg_{-1}=R(\varphi_1)$ yields $\fe(7)$.

\underline{{\bf Depth $d=2$}}. Here we need {\bf generalized}
prolongations, see \cite{Shch}. Again there are just a few algebras
$\fg_0$ and $\fg_0$-modules $\fg_{-1}$ for which $\fg_{1}\neq 0$ and
$\fg=\oplus \fg_{i}$ is simple:

\underline{$\fg(2)$; $\ff(4)$; $\fe(8)$}. These Lie algebras
correspond to the prolongations of their non-positive part (with
$\fg_0$ being isomorphic to $\fgl(2)$; $\fo(6)$ or $\fsp(6)$;
$\fe(7)$ or $\fo(14)$, respectively) in the following
$\Zee$-gradings. Let us rig the nodes of the Dynkin graph with
coefficients of linear dependence of the maximal root with respect
to simple ones. If any end node is rigged by a 2, mark it (just
one node even if several are rigged by 2's) and set the degrees of
the Chevalley generators to be: 
\begin{equation}\label{mark}
\deg X^\pm_i=\begin{cases}\pm 1&\text{if the $i$th node is marked}\\
0&\text{otherwise}.\end{cases} \end{equation}

\underline{$\fk(2n+1)$}. The cases where $(\fg_{-}, \fg_0)_*$ is
simple and of infinite dimension correspond to the prolongations of
the non-positive part of $\fsp(2n+2)$ in the $\Zee$-grading
(\ref{mark}) with the last node marked. Then $\fg_0=\fc\fsp(2n)$,
$\fg_{-1}=R(\varphi_1)$ and $\fg_{-1}$ is the trivial
$\fg_0$-module. In these cases, $(\fg_{-}, \fg_0)_*=\fk(2n+1)$.
\normalsize

\ssec{Superization}{}~{}

\sssec{Queerification} This is the functor  $Q: A\tto
Q(A):=A[\eps]$, where
$$ p(\eps)=\od, \; \eps^2=-1 \; \text{ and $ \eps a
=(-1)^{p(a)}a\eps$ for any $a\in A$}. $$ We set $\fq(n)=Q(\fgl(n))$.

\footnotesize

\sssec{Definition of Lie superalgebras for $p=2$} A Lie superalgebra
for $p=2$ as a superspace $\fg=\fg_\ev\oplus\fg_\od$ such that
$\fg_\ev$ is a Lie algebra, $\fg_\od$ is an $\fg_\ev$-module (made
into the two-sided one by symmetry; more exactly, by anti-symmetry,
but if $p=2$, it is the same) and on $\fg_\od$ a {\it
squaring}\index{squaring} (roughly speaking, the halved bracket) is
defined
\begin{equation}\label{squaring}
\begin{array}{c}
x\mapsto x^2\quad \text{such that $(ax)^2=a^2x^2$ for any $x\in
\fg_\od$ and $a\in \Kee$, and}\\
{}(x+y)^2-x^2-y^2\text{~is a bilinear form on $\fg_\od$ with values
in $\fg_\ev$.}
\end{array}
\end{equation}
For any $x, y\in \fg_\od$, we set
\begin{equation}\label{bracket} {}[x, y]:=(x+y)^2-x^2-y^2.
\end{equation}
We also assume, as usual, that

\medskip

if $x,y\in\fg_\ev$, then $[x,y]$ is the bracket on the Lie algebra;

if $x\in\fg_\ev$ and $y\in\fg_\od$, then
$[x,y]:=l_x(y)=-[y,x]=-r_x(y)$, where $l$ and $r$ are the left and
right $\fg_\ev$-actions on $\fg_\od$, respectively.
\medskip

The {\it Jacobi identity}\index{Jacobi identity} involving odd
element has now the following form:
\begin{equation}\label{JI}
~[x^2,y]=[x,[x,y]]\text{~for any~} x\in\fg_\od, y\in\fg.
\end{equation}

\normalsize

\begin{Conj}[Amended KL = Super KSh, $p>0$]
For $p>0$, to get all {\bf $\Zee$-graded} simple finite dimensional
examples of Lie algebras and Lie superalgebras:

(a) apply the KSh procedure to every simple Lie algebra of type (1)
over $\Cee$ (if $p=2$, apply the KSh procedure also to every simple
Lie superalgebra of of type (1) over $\Cee$ and their simple
Volichenko subalgebras described in \cite{LSer}),

(b) if $p=2$, apply queerification (as in \cite{Leb2}) to the
results of (a);

(c) if $p=2$, take Jurman's examples \cite{J} (and generalizations
of the same construction, if any: It looks like a specific $p=2$
non-super version of the queerification);

(d) take the non-positive part of every simple (up to center) finite
dimensional $\Zee$-graded algebra obtained at steps (a)--(c) and
(for $p=5$, $3$ and $2$) the exceptional ones of the form $\fg(A)$
listed in \cite{BGL5}, consider its complete and
partial\footnote{This term is too imprecise at the moment: it
embraces Frank and Ermolaev algebras, Frank and other exceptional
superalgebras (\cite{BGL6, BGL7}).} prolongs and distinguish their
simple subquotients.

To get non-graded examples, we have to take as a possible $\fg_0$
deformations of the simple algebras obtained at steps (a)--(d) and
Shen's \lq\lq variations'' \cite{Sh} (unless they can be interpreted
as deformations of the algebras obtained at earlier
steps).\end{Conj}

For preliminary results, see \cite{GL4, BjL},
\cite{BGL3}--\cite{BGL5}, \cite{ILL, Leb2}. (For $p=3$ and Lie
algebras, this is how Grozman and me got an interpretation of all
the \lq\lq mysterious" exceptional simple vectorial Lie algebras
known before \cite{GL4} was published; we also found two (if not
three) series of new simple algebras.) Having obtained a supply of
such examples, we can sit down to \so{compute certain cohomology in
order to describe their deformations} (provided we will be able to
understand what we are computing, cf. footnote 2); for the already
performed, see \cite{KKCh,KuCh,Ch, BGL4}.

\section{Further details}

\ssec{How to construct finite dimensional simple Lie algebras if
$p\geq 5$} Observe that although in the modular case there is a
wider variety of pairs $(\fg_{-1}, \fg_0)$ yielding nontrivial
prolongs than for $p=0$ (for the role of $\fg_0$ we can now take
vectorial Lie algebras or their central extensions), {\it a
posteriori} we know that we can always confine ourselves to the same
pairs $(\fg_{-1}, \fg_0)$ as for $p=0$. Melikyan's example looked as
a deviation from the pattern, but Kuznetsov's observation \cite{Ku1}
elaborated in \cite{GL4} shows that for $p\geq 5$ all is the same.
 Not so if
$p\leq 3$:

\ssec{New simple finite dimensional Lie algebras for $p=3$} In
\cite{S}, Strade listed known to him at that time examples of simple
finite dimensional Lie algebras for $p=3$. The construction of such
algebras is usually subdivided into the following types and
deformations of these types:

(1) algebras with Cartan matrix CM (sometimes encodable by Dynkin
graphs, cf. \cite{S, BGL5}),

(2) algebras of vectorial type (meaning that they have more roots
of one sign than of the other with respect to a partition into
positive and negative roots).

Case (1) was solved in \cite{WK, KWK}.

Conjecture 2 suggests to consider certain $\Zee$-graded prolongs
$\fg$. For Lie algebras and $p=3$, Kuznetsov described various
restrictions on the 0-th component of $\fg$ and the $\fg_0$-module
$\fg_{-1}$ (for partial summary, see \cite{GK, Ku1, Ku2}, \cite{BKK}
and a correction in \cite{GL4}). {\bf What are these restrictions
for  Lie algebras for $p=2$?  What are they for Lie superalgebras
for any $p>0$?}

\ssec{Exceptional simple finite dimensional Lie superalgebras for
$p\leq 5$}  Elduque investigated which spinor modules over
orthogonal algebras can serve as the odd part of a simple Lie
superalgebra and discovered an exceptional simple Lie superalgebra
for $p=5$. Elduque also superized the Freudenthal Magic Square and
expressed it in a new way, and his approach yielded nine new simple
(exceptional as we know now thanks to the classification
\cite{BGL5}) finite dimensional Lie superalgebras for $p=3$, cf.
\cite{CE, El1, CE2}. These Lie superalgebras possess Cartan matrices
(CM's) and we described all CMs and presentations of these algebras
in terms of Chevalley generators, see \cite{BGL5} and references in
it. In \cite{BGL5} 12 more examples of exceptional simple Lie
superalgebras are discovered; in \cite{BGL3}, we considered some of
their \lq\lq most promising" (in terms of prolongations)
$\Zee$-gradings and discovered several new series of simple
vectorial Lie superalgebras.

\ssec{New simple finite dimensional Lie algebras and Lie
superalgebras for $p=2$}  Lebedev \cite{Leb1, Leb2} offered a new
series of examples of simple orthogonal Lie algebras without CM.
Together with Iyer, we constructed their prolongations, missed in
\cite{Lin}, see \cite{ILL}; queerifications of these orthogonal
algebras are totally new examples of simple Lie superalgebras. CTS
prolongs of {\bf some}\footnote{We are unable to CTS the
superalgebras of dimension $>40$ on computers available to us,
whereas we need to be able to consider at least $250$.} of these
superalgebras and examples found in \cite{BGL5} are considered in
\cite{BGL7}.

\ssec{Conclusion}

Passing to Lie superalgebras we see that even their definition, as
well as that of their prolongations, are not quite straightforward
for $p=2$, but, having defined them (\cite{LL, Leb2}), it remains to
apply the above-described procedures to get at least a supply of
examples. To prove the completeness of the stock of examples for any
$p$ is a much more difficult task that requires serious preliminary
study of the representations of the examples known and to be
obtained --- more topics for Ph.D. theses.

The references currently in preparation are to be soon found in
\texttt{arXiv}.

\label{lastpage}


\begin{thebibliography}{999999}


\bibitem[BKK]{BKK}
Benkart, G.; Kostrikin, A. I.; Kuznetsov, M. I. The simple graded
Lie algebras of characteristic three with classical reductive
component $L\sb 0$. Comm. Algebra 24 (1996), no. 1, 223--234

\bibitem[BjL]{BjL}
Bouarroudj S., Leites D., Simple Lie superalgebras and
non-integrable distributions in characteristic $p$. Zapiski nauchnyh
seminarov POMI, t. 331 (2006), 15--29; Reprinted in J. Math. Sci.
(NY) 2007; \texttt{arXiv:math.RT/0606682}


\bibitem[BGL3]{BGL3}
Bouarroudj S., Grozman P., Leites D., New simple modular Lie
superalgebras as generalized prolongs. \texttt{arXiv:
math.RT/0704.0130}


\bibitem[BGL4]{BGL4}
Bouarroudj S., Grozman P., Leites D., Deformations of the somple
symmetric modular Lie superalgebras. \texttt{arXiv: }

\bibitem[BGL5]{BGL5}
Bouarroudj S., Grozman P., Leites D., Simple modular Lie
superalgebras with Cartan matrix. \texttt{arXiv: 0710.5149}

\bibitem[BGL6]{BGL6}
Bouarroudj S., Grozman P., Leites D., Simple modular Lie
superalgebras as partial prolongs. \texttt{arXiv: }

\bibitem[BGL7]{BGL7}
Bouarroudj S., Grozman P., Leites D., Leites D., New simple modular
Lie superalgebras as generalized prolongs in characteristic 2.
\texttt{arXiv: }

\bibitem[C]{C}
Cartan \'E., \"Uber die
einfachen Transformationsgrouppen, Leipziger Berichte (1893),
395--420. Reprinted in: {\em \OE uvres compl\`{e}tes}. Partie II.
(French) [Complete works. Part II] Alg\`{e}bre, syst\`{e}mes
diff\'erentiels et probl\`{e}mes d'\'equivalence. [Algebra,
differential systems and problems of equivalence] Second edition.
\'Editions du Centre National de la Recherche Scientifique (CNRS),
Paris, 1984.

\bibitem[Ch]{Ch}
Chebochko, N. G. Deformations of classical Lie algebras with a
homogeneous root system in characteristic two. I. (Russian) Mat.
Sb. 196 (2005), no. 9, 125--156; translation in Sb. Math. 196
(2005), no. 9-10, 1371--1402


\bibitem[Cla]{Cla}
Clarke B., Decomposition of the tensor product of two irreducible
$\mathfrak{sl}(2)$-modules in characteristic $3$,
MPIMiS\footnote{Available at \texttt{http://www.mis.mpg.de}.}
preprint 145/2006

\bibitem[CE]{CE}
Cunha I., Elduque A., An extended Freudenthal magic square in
characteristic $3$; \texttt{arXiv:math.RA/0605379}

\bibitem[CE2]{CE2}
Cunha I., Elduque, A., The extended Freudenthal Magic Square and
Jordan algebras; \texttt{arXiv:math.RA/0608191}


\bibitem[El1]{El1}
Elduque, A. New simple Lie superalgebras in characteristic 3. J.
Algebra 296 (2006), no. 1, 196--233

\bibitem[El2]{El2}
Elduque, A. Some new simple modular Lie superalgebras.
\texttt{arXiv:math.RA/0512654}

\bibitem[Er]{Er}
Ermolaev, Yu. B. Integral bases of classical Lie algebras.
(Russian) Izv. Vyssh. Uchebn. Zaved. Mat. 2004, , no. 3, 16--25;
translation in Russian Math. (Iz. VUZ) 48 (2004), no. 3, 13--22.


\bibitem[FH]{FH}
Fulton, W., Harris, J., {\em Representation theory. A first
course}. Graduate Texts in Mathematics, 129. Readings in
Mathematics. Springer-Verlag, New York, 1991. xvi+551 pp

\bibitem[GK]{GK}
Gregory, T.; Kuznetsov, M. On depth-three graded Lie algebras of
characteristic three with classical reductive null component.
Comm. Algebra 32 (2004), no. 9, 3339--3371

\bibitem[GG]{GG}
Grishkov, A.; Guerreiro, M. New simple Lie algebras over fields of
characteristic 2. Resenhas 6 (2004), no. 2-3, 215--221


\bibitem[Gr]{Gr} Grozman P., \texttt{SuperLie},
\texttt{http://www.equaonline.com/math/SuperLie}


\bibitem[GL1]{GL1}
Grozman P., Leites D., Defining relations for classical Lie
superalgebras with Cartan matrix,  Czech.  J. Phys., Vol. 51, 2001,
no.  1, 1--22; \texttt{arXiv: hep-th/9702073}

\bibitem[GL2]{GL2}
Grozman P., Leites D., \texttt{SuperLie} and problems (to be)
solved with it. Preprint MPIM-Bonn, 2003-39
(\texttt{http://www.mpim-bonn.mpg.de})



\bibitem[GL4]{GL4}
Grozman P., Leites D., Structures of $G(2)$ type and nonintegrable
distributions in characteristic $p$. Lett. Math. Phys.  74 (2005),
no. 3, 229--262; \texttt{arXiv:math.RT/0509400}

\bibitem[GLS]{GLS}
Grozman P., Leites D., Shchepochkina I., Invariant operators on
supermanifolds and standard models. In: In: M.~Olshanetsky,
A.~Vainstein (eds.)  {\em Multiple facets of quantization and
supersymmetry.  Michael Marinov Memorial Volume}, World Sci.
Publishing, River Edge, NJ, 2002, 508--555;
\texttt{arXiv:math.RT/0202193}

\bibitem[ILL]{ILL}
Iyer U., Lebedev A., Leites, D. Prolongations of orthogonal Lie
algebras in characteristic 2. IN PREPARATION



\bibitem[J]{J}
Jurman, G.,  A family of simple Lie algebras in characteristic
two.  J. Algebra 271 (2004), no. 2, 454--481.


\bibitem[K]{K}
Kac V. {\em Infinite-dimensional Lie algebras}. Third edition.
Cambridge University Press, Cambridge, 1990. xxii+400 pp.



\bibitem[K2]{K2}
Kac V., Lie superagebras, Adv.  Math.  v. 26, 1977, 8--96

\bibitem[K3]{K3}
Kac, V. Classification of supersymmetries.  Proceedings of the
International Congress of Mathematicians, Vol. I (Beijing, 2002),
Higher Ed. Press, Beijing, 2002, 319--344

Cheng, Shun-Jen; Kac, V., Addendum: ``Generalized Spencer
cohomology and filtered deformations of ${\Zee}$-graded Lie
superalgebras" [Adv. Theor. Math. Phys. 2 (1998), no. 5,
1141--1182; MR1688484 (2000d:17025)].  Adv. Theor. Math. Phys.  8
(2004), no. 4, 697--709.

Cantarini, N.; Cheng, S.-J.; Kac, V. Errata to: ``Structure of
some $\Zee$-graded Lie superalgebras of vector fields" [Transform.
Groups 4 (1999), no. 2-3, 219--272; MR1712863 (2001b:17037)] by
Cheng and Kac. Transform. Groups 9 (2004), no. 4, 399--400

\bibitem[KWK]{KWK}
Kac, V. G. Corrections to: "Exponentials in Lie algebras of
characteristic $p$" [Izv. Akad. Nauk SSSR 35 (1971), no. 4,
762--788; MR0306282 (46 \#5408)] by B. Yu. Veisfeiler and Kac.
(Russian) Izv. Ross. Akad. Nauk Ser. Mat. 58 (1994), no. 4, 224;
translation in Russian Acad. Sci. Izv. Math. 45 (1995), no. 1, 229


\bibitem[KKCh]{KKCh}
Kirillov, S. A.; Kuznetsov, M. I.; Chebochko, N. G. Deformations
of a Lie algebra of type $G\sb 2$ of characteristic three.
(Russian) Izv. Vyssh. Uchebn. Zaved. Mat. 2000, no. 3, 33--38;
translation in Russian Math. (Iz. VUZ) 44 (2000), no. 3, 31--36

\bibitem[KL]{KL}
Kochetkov Yu., Leites D., Simple finite dimensional Lie algebras
in characteristic 2 related to superalgebras and on a notion of
finite simple group.  In: L.~A.~Bokut, Yu.~L.~Ershov and
A.~I.~Kostrikin (eds.)  {\em Proceedings of the International
Conference on Algebra.  Part 1., Novosibirsk, August 1989},
Contemporary Math. 131, Part 1, 1992, AMS, 59--67

\bibitem[KSh]{KSh}
Kostrikin, A. I., Shafarevich,  I.R., Graded Lie algebras of
finite characteristic, Izv. Akad. Nauk. SSSR Ser. Mat. 33 (1969)
251--322 (in Russian); transl.: Math. USSR Izv. 3 (1969) 237--304



\bibitem[KuCh]{KuCh}
Kuznetsov, M. I.; Chebochko, N. G. Deformations of classical Lie
algebras. (Russian) Mat. Sb. 191 (2000), no. 8, 69--88;
translation in Sb. Math. 191 (2000), no. 7-8, 1171--1190


\bibitem[Ku1]{Ku1}
Kuznetsov, M. I. The Melikyan algebras as Lie algebras of the type
$G\sb 2$. Comm. Algebra 19 (1991), no. 4, 1281--1312.

\bibitem[Ku2]{Ku2}
Kuznetsov, M. I. Graded Lie algebras with the almost simple
component $L\sb 0$. Pontryagin Conference, 8, Algebra (Moscow,
1998). J. Math. Sci. (New York) 106 (2001), no. 4, 3187--3211.

\bibitem[Leb1]{Leb1}
Lebedev A., Non-degenerate bilinear forms in characteristic $2$,
related contact forms, simple Lie algebras and superalgebras.
\texttt{arXiv: math.AC/0601536}

\bibitem[Leb2]{Leb2}
Lebedev A., Simple modular Lie superalgebras. \texttt{arXiv:
math.RT/...to appear}


\bibitem[LL]{LL}
Lebedev A., Leites D., On realizations of the Steenrod algebras. J.
Prime Res. Math., v. 2 no. 1 (2006), 1--13; MPIMiS preprint 131/2006



\bibitem[LSer]{LSer}
Leites, D.; Serganova, V. Metasymmetry and Volichenko algebras.
Phys. Lett. B 252 (1990), no. 1, 91--96; id, Symmetries wider than
supersymmetry. In: Dupllij S., Wess J. (eds.) {\em Noncommutative
structures in mathematics and physics} (Kiev, 2000), 13--30, NATO
Sci. Ser. II Math. Phys. Chem., 22, Kluwer Acad. Publ., Dordrecht,
2001

\bibitem[LSh]{LSh}
Leites D., Shchepochkina I., Classification of the simple Lie
superalgebras of vector fields, preprint\footnote{Available at
\texttt{http://www.mpi-bonn.mpg.de}.}  MPIM-2003-28


\bibitem[Lvd]{Lvd}
Leur J. van de, {\em Contragredient Lie superalgebras of finite
growth} (Ph.D. thesis) Utrecht, 1986; a short version published in
Commun.  in Alg., v.  17, 1989, 1815--1841

\bibitem[Lin]{Lin}
Lin, L., Nonalternating hamiltonian algebra $P(n, m)$ of
characteristic two, Comm. Algebra 21 (2) (1993) 399--411

\bibitem[Lin1]{Lin1}
Lin, L., {\em Lie algebras $K({\cF}, \mu_i)$ of Cartan type of
characteristic $p=2$ and their subalgebras} (Chinese. English
summary), J. East China Norm. Univ. Natur. Sci. Ed. 1 (1988), 16--23
MR0966993 (89k:17033)


\bibitem[RSh]{RSh}
Rudakov, A. N.; Shafarevich, I. R. Irreducible representations of
a simple three-dimensional Lie algebra over a field of finite
characteristic. (Russian) Mat. Zametki 2 (1967) 439--454 (1967, 2,
760--767)

\bibitem[Se]{Se}
Serganova, V., Automorphisms of simple Lie superalgebras.
(Russian) Izv. Akad. Nauk SSSR Ser. Mat. 48 (1984), no. 3,
585--598

\bibitem[Se1]{Se1}
Serganova, V., On generalizations of root systems. Comm. Algebra
24 (1996), no. 13, 4281--4299

\bibitem[Ser]{Ser}
Sergeev A., Orthogonal polynomials and Lie superalgebras,
\texttt{arXiv:math.RT/9810110}.


\bibitem[Shch]{Shch}
Shchepochkina I., How to realize Lie algebras by vector fields.
Theor. Mat. Fiz. 147 (2006) no. 3, 821--838;
\texttt{arXiv:math.RT/0509472}

\bibitem[Sh14]{Sh14}
Shchepochkina I., Five exceptional simple Lie superalgebras of
vector fields and their fourteen regradings. Representation Theory
(electronic journal of AMS), v. 3, 1999, 3 (1999), 373--415;
\texttt{arXiv:hep-th/9702121}

\bibitem[Sh]{Sh}
Shen, Guang Yu, Variations of the classical Lie algebra $G\sb 2$
in low characteristics. Nova J. Algebra Geom. 2 (1993), no. 3,
217--243.

\bibitem[Sk]{Sk}
Skryabin, S. M. New series of simple Lie
algebras of characteristic $3$. (Russian. Russian summary) Mat.
Sb. 183 (1992), no. 8, 3--22; translation in Russian Acad. Sci.
Sb. Math. 76 (1993), no. 2, 389--406

\bibitem[St]{St}
Steinberg, R. {\em Lectures on Chevalley groups}. Notes prepared
by John Faulkner and Robert Wilson. Yale University, New Haven,
Conn., 1968. iii+277 pp.


\bibitem[S]{S}
Strade, H. {\em Simple
Lie algebras over fields of positive characteristic. I. Structure
theory.} de Gruyter Expositions in Mathematics, 38. Walter de
Gruyter \& Co., Berlin, 2004. viii+540 pp.

\bibitem[WK]{WK}
Weisfeiler, B. Ju.; Kac, V. G. Exponentials in Lie algebras of
characteristic $p$. (Russian) Izv. Akad. Nauk SSSR Ser. Mat. 35
(1971), 762--788.


\bibitem[Y]{Y}
Yamaguchi K., Differential systems associated with simple graded
Lie algebras. Progress in differential geometry, Adv. Stud. Pure
Math., 22, Math. Soc. Japan, Tokyo, 1993, 413--494

\end{thebibliography}
\end{document}